\documentclass[11pt]{amsart}
\usepackage{amssymb}
\usepackage{mathrsfs}
\usepackage{enumerate}
\makeatletter
\@namedef{subjclassname@2010}{%
  \textup{2010} Mathematics Subject Classification}
\makeatother
\newtheorem{lemma}{Lemma}[section]
\newtheorem{theorem}{Theorem}[section]
\newtheorem{proposition}{Proposition}[section]
\newtheorem{corollary}{Corollary}[section]
\newtheorem{definition}{Definition}[section]

\newtheorem{example}{Example}[section]

\begin{document}

\title{ON CONTINUOUS WEAVING FUSION FRAMES IN HILBERT SPACES}

\author[V. Sadri]{Vahid Sadri}
\address{Institute of Fundamental Sciences\\University of Tabriz\\, Iran\\}
\email{vahidsadri57@gmail.com}

\author[R. Ahmadi]{Reza Ahmadi}%{Corresponding Author.}
\address{Institute of Fundamental Sciences\\University of Tabriz\\, Iran\\}
\email{rahmadi@tabrizu.ac.ir}

\author[Gh. Rahimlou]{Gholamreza Rahimlou}
\address{Department of Mathematics, Faculty of Tabriz  Branch\\ Technical and Vocational University (TUV), East Azarbaijan
, Iran\\}
\email{grahimlou@gmail.com}
%\date{}
\begin{abstract}
In this note, we first introduce the notation of weaving c-fusion frames in separable Hilbert spaces. After reviewing the conditions for maintaining the weaving c-fusion frames under the bounded linear operator and also, removing vectors from these frames, we will present a necessarily and sufficient condition about c-woven and c-fusion woven. Finally, perturbation of these frames will be introduced.
\end{abstract}

\subjclass[2010]{Primary 42C15; Secondary 42C40,41A58}

\keywords{Fusion frame, fusion woven, continuous weaving frame.}

\maketitle

\section{Introduction and Perliminaries}
Nowadays, frames (or discrete frames) have a significant role in both pure and applied mathematics, so that these are a fundamental research area in mathematics, computer science and engineering. Frames were introduced by Duffin and Scheaffer \cite{ds} in the context of non-harmonic Fourier series and since then, there have been many generalizations such as c-frame, g-frame, fusion frame, $K$-frame and etc.

 Continuous frames (or briefly c-frames) were proposed by Kaiser \cite{kaiser} and also independently by Ali et al. \cite{ali} to a family indexed by some locally compact space endowed with a Radon measure. C-frames are the first  generalizations frames to measure spaces. For more studies about these frames, we refer to \cite{fa3, rah2}.

Casazza and Kutyniok \cite{ck} were able to introduce fusion frames (or frame of subspaces) which are an important generalization of frames. Recent studies shows that fusion frames provide effective framework for modeling of sensor network, signal and image processing, sampling theory, filter bank and a variety of application that cannot  be modeled by discrete frames. After, Faroughi et al. \cite{fa1} introduced c-fusion frames which have been obtained from the combination of fusion and continuous frames. For more details, we refer to \cite{ah, fa, nr, rah}.

Recently, Bemrose et al. \cite{bem} introduced a new concept of weaving frames which is motivated by a question in distributed signal processing. In \cite{arab, vash2}, these frames have been presented for fusion frame and also, Vashisht and Deepshikha \cite{vash} were able to introduce for continuous case. In this paper, we will present weaving for c-fusion frames.

Throughout this paper, $(X,\mu)$ is a measure space with positive measure $\mu$, $H$ is a Hilbert space, $\mathbb{H}$ is the collection of all closed subspaces of $H$, $\pi_{V}$ is the orthogonal projection from $H$ onto a closed subspace $V$ and  $\mathcal{B}(H,K)$ is the set of all bounded and linear operators from $H$ to $K$. If $H=K$, then $\mathcal{B}(H,H)$ will be denoted by $\mathcal{B}(H)$. For each $m>1$, we define $[m]:=\{1,2,\cdots, m\}$.

We present some theorems in operator theory which will be needed in the next sections.
\begin{lemma}(\cite{ga})\label{l1}
Let $V\subseteq H$ be a closed subspace, and $U$ be a linear  bounded operator on $H$. Then
\begin{equation*}
\pi_{V}U^*=\pi_{V}U^* \pi_{\overline{UV}}.
\end{equation*}
If $U$ is an unitary (i.e. $U$ is bijective and $U^*=U^{-1}$), then
$\pi_{\overline{UV}}U=U\pi_{V}$.
\end{lemma}
\begin{lemma}[\cite{dag}]\label{dag}
Let $L_1\in\mathcal{B}(H_1,H)$ and $L_2\in\mathcal{B}(H_2,H)$. Then the following assertions are equivalent:
\begin{enumerate}
\item[(I)] $\mathcal{R}(L_1)\subseteq\mathcal{R}(L_2)$;
\item[(II)] $L_1L_1^*\leq\lambda L_2L_2^*$ for some$\lambda>0$;
\item[(III)] there exists  $U\in\mathcal{B}(H_1,H_2)$ such that $L_1=L_2 U$.
\end{enumerate}

Moreover, if those conditions are valid then there exists a unique operator $U$ such that
\begin{enumerate}
\item[(a)] $\left\| U\right\|^2=\inf\{\alpha>0 \ \vert \ L_1L_1^*\leq\alpha L_2L_2^*\};$
\item[(b)] $\ker{L_1}=\ker{U};$
\item[(c)] $\mathcal{R}(U)\subseteq\overline{\mathcal{R}(L_2^*)}.$
\end{enumerate}
\end{lemma}
If an operator $U$ has closed range, then there exists a right-inverse operator $U^ \dagger$ (pseudo-inverse of $U$) in the following sense (see \cite{ch}):
\begin{lemma}\label{dager}
Let $U\in\mathcal{B}(K,H)$  be a bounded operator with closed range $\mathcal{R}(U)$. Then there exists a bounded operator $U^\dagger \in\mathcal{B}(H,K)$ for which
\begin{equation*}
UU^{\dagger} x=x, \ \ x\in \mathcal{R}(U).
\end{equation*}
\end{lemma}
\begin{lemma}\label{Ru}
Let $U\in\mathcal{B}(K,H)$. Then the following assertions holds:
\begin{enumerate}
\item[(I)] $\mathcal{R}(U)$ is closed in $H$ if and only if $\mathcal{R}(U^{\ast})$ is closed in $K$;
\item[(II)] $(U^{\ast})^\dagger=(U^\dagger)^\ast$;
\item[(III)] The orthogonal projection of $H$ onto $\mathcal{R}(U)$ is given by $UU^{\dagger}$;
\item[(IV)] The orthogonal projection of $K$ onto $\mathcal{R}(U^{\dagger})$ is given by $U^{\dagger}U$;
\item[(V)] On $\mathcal{R}(U)$, the operator $U^{\dagger}$ is given explicitly by $U^{\dagger}=U^*(U^*)^{-1}$.
\end{enumerate}
\end{lemma}

\begin{definition}[c-frame]
Let $F:X\rightarrow H$ be a mapping such that the mapping $x\rightarrow\langle h,F(x)\rangle$ of $X$ to $\Bbb C$ is measurable (i.e. weakly measurable) for each $h\in H$. $F$ is called c-frame for $H$ (with respect to $\mu$) if there exists $0<A\leq B<\infty$ such that for every $h\in H$
\begin{eqnarray*}
A\Vert h\Vert^2\leq\int_{X}\vert\langle h,F(x)\rangle\vert^2\,d\mu\leq B\Vert h\Vert^2.
\end{eqnarray*}
\end{definition}
 Suppose that $F:X\rightarrow\mathbb{H}$ and we denote $\mathscr{L}^2(X,F)$ the class of all weakly measurable mappings $f:X\rightarrow H$ (i.e. for all $h\in H$, the mapping $x\rightarrow\left\langle f(x), h\right\rangle$ is measurable) such that for any $x\in X$, $f(x)\in F(x)$ and
\begin{equation*}
\int_{X}\left\| f(x)\right\|^2\,d\mu<\infty.
\end{equation*}
It can be a Hilbert space with the inner product defined by
\begin{equation*}
\left\langle f,g\right\rangle=\int_{X}\left\langle f(x),g(x)\right\rangle\,d\mu, \quad f,g\in\mathscr{L}^2(X,F).
\end{equation*}
\begin{definition}[c-fusion frame]
Assume that $F:X\rightarrow\mathbb{H}$  such that for each $h\in H$, the mapping $x\rightarrow\pi_{F(x)}(h)$ is measurable (i.e. is weakly measurable) and $v:X\rightarrow\Bbb R^+$ be a measurable function. Then $(F,v)$ is called a continuous fusion frame (or c-fusion frame)  for $H$ (with respect to $\mu$) if
there exist  $0<A\leq B<\infty$ such that for each $h\in H$
\begin{equation}\label{f1}
A\Vert h\Vert^2\leq\int_{X}v^2(x)\Vert \pi_{F(x)}h\Vert^2\,d\mu\leq B\Vert h\Vert^2.
\end{equation}
\end{definition}
When the right hand side of (\ref{f1}) holds, $(F,v)$ is called a c-fusion Bessel sequence for $H$ with  bound $B$.
We say $(F,v)$ is a Parseval c-fusion frame whenever $A=B=1$.
The synthesis operator is defined  weakly as follows (for more details we refer \cite{fa1}):
\begin{align*}
T_{F}&:\mathscr{L}^2(X,F) \to H,\\
T_F&(f)=\int_{X}vf\,d\mu,
\end{align*}
where for each $h\in H$
\begin{equation*}
\left\langle T_{F}(f),h\right\rangle=\int_{X} v(x)\left\langle f(x),h\right\rangle\,d\mu.
\end{equation*}
The analysis operator is given by
 \begin{align*}
T^{*}_{F}:H&\longrightarrow\mathscr{L}^2(X,F),\\
T^{*}_{F}(h)&=v\pi_{F}h.
\end{align*}
Finally, the  c-fusion frame operator $S_{F}:=T_{F}T^*_{F}$ is defined by
\begin{align*}
S_{F}&:H\longrightarrow H,\\
S_F(h)&=\int_{X}v^2\pi_{F}(h)\,d\mu.
\end{align*}
Hence for each $h_1, h_2 \in H$
\begin{equation*}
\left\langle S_F(h_1),h_2\right\rangle=\int_{X} v^2(x)\left\langle\pi_{F(x)}h_1,h_2 \right\rangle\,d\mu.
\end{equation*}
Therefore,
\begin{equation*}
AId_H\leq S_F\leq BId_H
\end{equation*}
and we obtain, if $(F,v)$ is a c-fusion frame, then $S_F$ is a positive, self-adjoint and invertible operator.
\begin{definition}[weaving frame]
A family of  frames $\{f_{ij}\}_{j\in\Bbb J, i\in [m]}$ for $H$ which $\Bbb J\subseteq\Bbb Z$,  is said to be  woven if there exist universal same positive constants $0<A\leq B$ such that for each partition $\{\sigma_i\}_{i\in [m]}$ of $\Bbb J$, the family $\{f_{ij}\}_{j\in \sigma_i, i\in [m]}$ is a  frame for $H$ with bounds $A$ and $B$. Each family $\{f_{ij}\}_{j\in\Bbb J, i\in [m]}$ is called a \textit{weaving}.
\end{definition}
\begin{example}(\cite{bem})
Let $\varepsilon>0$, set $\delta=(1+\varepsilon^2)^{-\frac{1}{2}}$ and let $\{e_1, e_2\}$ be the standard orthonormal basis of $\Bbb R^2$. Then the two sets
$$\Phi=\{\varphi_i\}_{i=1}^{4}=\{\delta e_1, \delta\varepsilon e_1, \delta e_2, \delta\varepsilon e_2\}$$
and
$$\Psi=\{\psi_i\}_{i=1}^{4}=\{\delta\varepsilon e_1, \delta e_1, \delta\varepsilon e_2, \delta e_2,\}$$
are Parseval frames, which are woven since each choice of $\sigma$ gives a spanning set.
\end{example}
%%%%%%%%%%%%%%%%%%%%%%%%%%%%%%%%%
\section{Continuous Weaving Fusion Frames}
Throughout the paper, by partition of a measure space $(X,\mu)$ we mean partition of $X$ into disjoint measurable
sets.
\begin{definition}
A family of c-fusion frames $(F_i, v_i)_{i\in[m]}$  for $H$ is said to be continuous fusion woven (or c-fusion woven) if there exist universal positive constants $0<A\leq B$ such that for each partition $\{\sigma_i\}_{i\in [m]}$ of $X$, the family $\big(F_i(x),v_i(x)\big)_{i\in[m], x\in\sigma_i}$ is a c-fusion frame  for $H$ with bounds $A$ and $B$.
\end{definition}
In above definition, $A$ and $B$ is called universal c-fusion frame bounds.
The following Proposition shows that every c-fusion woven has an universal upper c-fusion frame bound.
\begin{proposition}\label{t1}
Let $(F_i, v_i)_{i\in[m]}$ be a c-fusion Bessel sequence for $H$ with bound $B_i$ for each $i\in [m]$. Then, for any partition $\{\sigma_i\}_{i\in [m]}$ of $X$, the family $\big(F_i(x),v_i(x)\big)_{i\in[m], x\in\sigma_i}$ is a c-fusion Bessel sequence with the Bessel bound $\sum_{i\in [m]}B_i$.
\end{proposition}
\begin{proof}
Let $\{\sigma_i\}_{i\in [m]}$ be any partition of $X$. For each $f\in H$, we have
\begin{align*}
\sum_{i\in[m]}\int_{\sigma_i}v_i^2(x)\Vert \pi_{F_i(x)}h\Vert^2\,d\mu\leq\sum_{i\in[m]}\int_{X}v_i^2(x)\Vert \pi_{F_i(x)}h\Vert^2\,d\mu\leq\Big(\sum_{i\in[m]}B_i\Big)\Vert f\Vert^2.
\end{align*}
\end{proof}
In next results, we construct a c-fusion woven by using a bounded linear operator.
\begin{theorem}
Let $\big(F_i(x),v_i(x)\big)_{i\in[m], x\in\sigma_i}$ be a c-fusion woven for $H$ with universal bounds $A, B$. If $U\in\mathcal{B}(H)$ has closed range, then $\big(\overline{UF_i}(x),v_i(x)\big)_{i\in[m], x\in\sigma_i}$ is a c-fusion woven for $\mathcal{R}(U)$ with frame bounds
$$A\Vert U^{\dagger}\Vert^{-2}\Vert U\Vert^{-2},\quad B\Vert U^{\dagger}\Vert^{2}\Vert U\Vert^2.$$
\end{theorem}
\begin{proof}
Let $h\in\mathcal{R}(U)$ and $U\in\mathcal{B}(H)$ be closed range. By
$$\pi_{F_i(x)}U^*=\pi_{F_i(x)}U^*\pi_{\overline{UF_i(x)}},$$
for all $x\in\sigma_i$ and $i\in[m]$,  the mapping $x\mapsto\pi_{\overline{UF_i(x)}}$ is weakly measurable. We have by Lemma \ref{l1}
\begin{align*}
A\Vert h\Vert^2&=A\Vert(U^{\dagger})^*U^*h\Vert^2\\
&\leq A\Vert U^{\dagger}\Vert^2\Vert U^*h\Vert^2\\
&\leq\Vert U^{\dagger}\Vert^2\sum_{i\in[m]}\int_{X}v_i^2(x)\Vert\pi_{F_i(x)}U^*h\Vert^2\,d\mu\\
&=\Vert U^{\dagger}\Vert^2\sum_{i\in[m]}\int_{X}v_i^2(x)\Vert\pi_{F_i(x)}U^*\pi_{\overline{UF_i(x)}}h\Vert^2\,d\mu\\
&\leq\Vert U^{\dagger}\Vert^2\Vert U\Vert^2\sum_{i\in[m]}\int_{X}v_i^2(x)\Vert\pi_{\overline{UF_i(x)}}h\Vert^2\,d\mu
\end{align*}
and the lower bound is evident. For the upper bound, we can write for all $h\in\mathcal{R}(U)$,
\begin{align*}
\sum_{i\in[m]}\int_{X}v_i^2(x)\Vert\pi_{\overline{UF_i(x)}}h\Vert^2\,d\mu&=\sum_{i\in[m]}\int_{X}v_i^2(x)\Vert\pi_{\overline{UF_i(x)}}(U^{\dagger})^*U^*h\Vert^2\,d\mu\\
&=\sum_{i\in[m]}\int_{X}v_i^2(x)\Vert\pi_{\overline{UF_i(x)}}(U^{\dagger})^*\pi_{U^{\dagger}\overline{UF_i(x)}}U^*h\Vert^2\,d\mu\\
&\leq\Vert U^{\dagger}\Vert^2\sum_{i\in[m]}\int_{X}v_i^2(x)\Vert\pi_{F_i(x)}U^*h\Vert^2\,d\mu\\
&\leq B\Vert U^{\dagger}\Vert^{2}\Vert U\Vert^2\Vert h\Vert^2.
\end{align*}
\end{proof}
\begin{corollary}
Let $\big(F_i(x),v_i(x)\big)_{i\in[m], x\in\sigma_i}$ be a c-fusion woven for $H$ with  bounds $A, B$. If $U\in\mathcal{B}(H)$ is an invertible operator, then $\big(UF_i(x),v_i(x)\big)_{i\in[m], x\in\sigma_i}$ is a c-fusion woven for $\mathcal{R}(U)$ with frame bounds
$$A\Vert U^{-1}\Vert^{-2}\Vert U\Vert^{-2},\quad B\Vert U^{-1}\Vert^{2}\Vert U\Vert^2$$
\end{corollary}
\begin{theorem}
Let $(F,v)$ and $(G,w)$ be c-fusion woven for $H$ and $W\subset H$ be a closed subspace. Then $(F(x)\cap W,v)$ and $(G(x)\cap W,w)$ are c-fusion woven  for $W$ for all $x\in X$.
\end{theorem}
\begin{proof}
Suppose that $h\in W$ and $\sigma\in X$ is a measurable subset. Then,
\begin{align*}
&\int_{\sigma}v^2(x)\Vert\pi_{F(x)}h\Vert^2\,d\mu(x)+\int_{\sigma^c}w^2(x)\Vert\pi_{G(x)}h\Vert^2\,d\mu(x)\\
&\quad=\int_{\sigma}v^2(x)\Vert\pi_{F(x)}\pi_{W}h\Vert^2\,d\mu(x)+\int_{\sigma^c}w^2(x)\Vert\pi_{G(x)}\pi_{W}h\Vert^2\,d\mu(x)\\
&\quad=\int_{\sigma}v^2(x)\Vert\pi_{F(x)\cap W}h\Vert^2\,d\mu(x)+\int_{\sigma^c}w^2(x)\Vert\pi_{G(x)\cap W}h\Vert^2\,d\mu(x).
\end{align*}
\end{proof}

The next proposition shows that it is enough to check c-weaving fusion on smaller measurable space than the original which this is an extension of Proposition 3.10 in \cite{vash}.
\begin{theorem}
For all $i\in[m]$, let  $(F_i(x), v_i(x))_{x\in X}$ be a c-fusion frame for $H$ with frame bounds $A_i$ and $B_i$. If there exists a measurable subset $Y\subset X$  such that the family of c-fusion frame $(F_i(x), v_i(x))_{x\in Y, i\in[m]}$ is a c-fusion woven for $H$ with universal frame bounds $A$ and $B$. Then, $(F_i(x), v_i(x))_{x\in X, i\in[m]}$ is a c-fusion woven for $H$ with universal bounds $A$ and $\sum_{i\in[m]}B_i$.
\end{theorem}
\begin{proof}
Suppose that $\{\sigma_i\}_{i\in[m]}$ is a partition of $X$. We define for each $h\in H$,
\begin{align*}
\varphi&:X\longrightarrow \Bbb C,\\
\varphi(x)&=\sum_{i\in[m]}\chi_{\sigma_i}(x)\Vert\pi_{F_i(x)}h\Vert.
\end{align*}
Therefore, $\varphi$ is measurable. For any $h\in H$, we have
\begin{align*}
\sum_{i\in[m]}\int_{\sigma_i}v_i^2(x)\Vert\pi_{F_i(x)}h\Vert^2\,d\mu(x)
&\leq\sum_{i\in[m]}\int_{X}v_i^2(x)\Vert\pi_{F_i(x)}h\Vert^2\,d\mu(x)\\
&\leq\big(\sum_{i\in[m]} B_i\big)\Vert h\Vert^2.
\end{align*}
For the lower bound, it is clear that $\{\sigma_i\cap Y\}_{i\in[m]}$ is a partition of $Y$. Thus, $(F_i(x), v_i(x))_{x\in\sigma_i\cap Y, i\in[m]}$ is a c-fusion frame for $H$ with the lower frame bound $A$. Hence,
\begin{align*}
\sum_{i\in[m]}\int_{\sigma_i}v_i^2(x)\Vert\pi_{F_i(x)}h\Vert^2\,d\mu(x)
&\geq\sum_{i\in[m]}\int_{\sigma_i\cap Y}v_i^2(x)\Vert\pi_{F_i(x)}h\Vert^2\,d\mu(x)\\
&\geq A\Vert h\Vert^2.
\end{align*}
\end{proof}
Casazza and Lynch in \cite{cly} showed that It is possible to remove vectors from woven frames and still
be left with woven frames. After, this topic was presented in \cite{vash} and now, we study it in the following Theorem.
\begin{theorem}
Let  $(F_i(x), v_i(x))_{x\in X, i\in[m]}$ be a c-fusion woven for $H$ with universal bounds $A$ and $B$. If there exists $0<D<A$ and a measurable subset $Y\subset X$ and $n\in[m]$ such that
$$\sum_{i\in[m]\setminus\{n\}}\int_{X\setminus Y}v_i^2(x)\Vert\pi_{F_i(x)}h\Vert^2\,d\mu(x)\leq D\Vert h\Vert^2  \qquad \forall h\in H,$$
then, the family $(F_i(x), v_i(x))_{x\in Y, i\in[m]}$ is a c-fusion woven for $H$ with frame bounds $A-D$ and $B$.
\end{theorem}
\begin{proof}
Suppose that $\{\sigma_i\}_{i\in[m]}$ is a partition of $Y$ and $\{\tau_i\}_{i\in[m]}$ is a partition of $X\setminus Y$. For a given $h\in H$, we define
\begin{align*}
\varphi&:Y\longrightarrow \Bbb C,\\
\varphi(x)&=\sum_{i\in[m]}\chi_{\sigma_i}(x)\Vert\pi_{F_i(x)}h\Vert,
\end{align*}
and
\begin{align*}
\phi&:X\longrightarrow \Bbb C,\\
\phi(x)&=\sum_{i\in[m]}\chi_{\sigma_i\cup\tau_i}(x)\Vert\pi_{F_i(x)}h\Vert.
\end{align*}
Since $(F_i(x), v_i(x))_{x\in\sigma_i\cup\tau_i, i\in[m]}$ is a c-fusion frame for $H$ and $\varphi=\phi|_{Y}$, then $\varphi$ and $\phi$ are measurable. So, for each $h\in H$, we have
\begin{align*}
\sum_{i\in[m]}\int_{\sigma_i}v_i^2(x)\Vert\pi_{F_i(x)}h\Vert^2\,d\mu(x)
&\leq\sum_{i\in[m]}\int_{\sigma_i\cup\tau_i}v_i^2(x)\Vert\pi_{F_i(x)}h\Vert^2\,d\mu(x)\\
&\leq B\Vert h\Vert^2.
\end{align*}
Now, for the lower frame bound, assume that $\{\varsigma_i\}_{i\in[m]}$ such that $\varsigma_n=\emptyset$. Then, $\{\varsigma_i\cup\sigma_i\}_{i\in[m]}$ is a partition of $X$ and so, for any $h\in H$,
\begin{align*}
&\sum_{i\in[m]}\int_{\sigma_i}v_i^2(x)\Vert\pi_{F_i(x)}h\Vert^2\,d\mu(x)\\
&\quad=\sum_{i\in[m]\setminus\{n\}}\bigg(\int_{\varsigma_i\cup\sigma_i}v_i^2(x)\Vert\pi_{F_i(x)}h\Vert^2\,d\mu(x)
-\int_{\varsigma_i}v_i^2(x)\Vert\pi_{F_i(x)}h\Vert^2\,d\mu(x)\bigg)\\
&\qquad +\int_{\sigma_n}v_i^2(x)\Vert\pi_{F_i(x)}h\Vert^2\,d\mu(x)\\
&\quad\geq\sum_{i\in[m]\setminus\{n\}}\bigg(\int_{\varsigma_i\cup\sigma_i}v_i^2(x)\Vert\pi_{F_i(x)}h\Vert^2\,d\mu(x)
-\int_{X\setminus Y}v_i^2(x)\Vert\pi_{F_i(x)}h\Vert^2\,d\mu(x)\bigg)\\
&\qquad +\int_{\sigma_n}v_i^2(x)\Vert\pi_{F_i(x)}h\Vert^2\,d\mu(x)\\
&\quad=\sum_{i\in[m]}\int_{\varsigma_i\cup\sigma_i}v_i^2(x)\Vert\pi_{F_i(x)}h\Vert^2\,d\mu(x)-
\sum_{i\in[m]\setminus\{n\}}\int_{X\setminus Y}v_i^2(x)\Vert\pi_{F_i(x)}h\Vert^2\,d\mu(x)\\
&\quad\geq(A-D)\Vert h\Vert^2.
\end{align*}
\end{proof}

\begin{theorem}\label{tfF}
Let $v_i,w_i:X\rightarrow\Bbb R^+$ be measurable functions, $X_i\subseteq X$ be measurable subset for all $i\in[m]$ and $\{F_{i}\}_{i\in[m]}$, $\{G_{i}\}_{i\in[m]}$ be c-frame sequences on $X_i$ for $H$ with frame bounds $(A_{F_i}, B_{F_i})$ and $(A_{G_i}, B_{G_i})$, respectively. Assume that $\mathcal{F}_i,\mathcal{G}_i:X\rightarrow\mathbb{H}$ for any $i\in[m]$ such that
\begin{align*}
\mathcal{F}_i(x):=\begin{cases}
\overline{\mathrm{span}}\{F_i(x)\},& x\in X_i\\
0,& x\notin X_i
\end{cases}, \qquad \mathcal{G}_i(x):=\begin{cases}
\overline{\mathrm{span}}\{G_i(x)\},& x\in X_i\\
0,& x\notin X_i
\end{cases}
\end{align*}
and also
\begin{align*}
&0<A_F:=\inf_{i\in\Bbb I}A_{F_i}\leq B_F:=\sup_{i\in\Bbb I}B_{F_i}<\infty,\\ &0<A_G:=\inf_{i\in\Bbb I}A_{G_i}\leq B_G:=\sup_{i\in\Bbb I}B_{G_i}<\infty.
\end{align*}
Let for all $i\in[m]$,
\begin{align*}
\mathfrak{F}_i, \mathfrak{G}_i&: X\times X_i\longmapsto H,\\
\mathfrak{F}_i(x,y)&=F_i(y),\quad \mathfrak{G}_i(x,y)=G_i(y),
\end{align*}
and
\begin{align*}
\mathfrak{v}_i, \mathfrak{w}_i&: X\times X_i\longmapsto \Bbb R^+,\\
\mathfrak{v}_i(x,y)&=v_i(x),\quad \mathfrak{w}_i(x,y)=w_i(x).
\end{align*}
Then the following assertions are equivalent.
\begin{enumerate}
\item[(I)] $\{\mathfrak{F}_i.\mathfrak{v}_i\}_{i\in[m]}$ and $\{\mathfrak{G}_i.\mathfrak{w}_i\}_{i\in[m]}$ are c-woven for $H$.
\item[(II)] $(\mathcal{F}_i, v_i)_{i\in[m]}$ and $(\mathcal{G}_i, w_i)_{i\in[m]}$ are c-fusion woven for $H$.
\end{enumerate}
\end{theorem}
\begin{proof}
$(I)\Rightarrow(II)$. Suppose that $\sigma\subset X$ is a measurable subset and $h\in H$. Let $A:=\min\{A_{F}, A_{G}\}$ and $\{\mathfrak{F}_i.\mathfrak{v}_i\}_{i\in[m]}$,$\{\mathfrak{G}_i.\mathfrak{w}_i\}_{i\in[m]}$ are c-woven for $H$ with universal frame bounds $C, D$ for each $x\in X$ and $y\in X_i$. We have
\begin{align*}
&A\int_{\sigma}v_i^2(x)\Vert\pi_{\mathcal{F}_i(x)}h\Vert^2\,d\mu(x)+A\int_{\sigma^c}w_i^2(x)\Vert\pi_{\mathcal{G}_i(x)}h\Vert^2\,d\mu(x)\\
&\quad \leq \int_{\sigma}A_{F_i}\Vert v_i(x)\pi_{\mathcal{F}_i(x)}h\Vert^2\,d\mu(x)+\int_{\sigma^c}A_{G_i}\Vert w_i(x)\pi_{\mathcal{G}_i(x)}h\Vert^2\,d\mu(x)\\
&\quad \leq \int_{\sigma}\int_{X_i}\vert\langle v_i(x)\pi_{\mathcal{F}_i(x)}h, F_i(y)\rangle\vert^2\,d\mu(y)\,d\mu(x)\\
&\qquad +\int_{\sigma^c}\int_{X_i}\vert\langle w_i(x)\pi_{\mathcal{G}_i(x)}h, G_i(y)\rangle\vert^2\,d\mu(y)\,d\mu(x)\\
&\quad =\int_{\sigma}\int_{X_i}\vert\langle \pi_{\mathcal{F}_i(x)}h, v_i(x)F_i(y)\rangle\vert^2\,d\mu(y)\,d\mu(x)\\
&\qquad +\int_{\sigma^c}\int_{X_i}\vert\langle \pi_{\mathcal{G}_i(x)}h, w_i(x)G_i(y)\rangle\vert^2\,d\mu(y)\,d\mu(x)\\
&\quad =\int_{\sigma}\int_{X_i}\vert\langle h, v_i(x)F_i(y)\rangle\vert^2\,d\mu(y)\,d\mu(x)\\
&\qquad +\int_{\sigma^c}\int_{X_i}\vert\langle h, w_i(x)G_i(y)\rangle\vert^2\,d\mu(y)\,d\mu(x)\\
&\quad\leq D\Vert h\Vert^2.
\end{align*}
With the same way, we conclude that
\begin{align*}
&B\int_{\sigma}v_i^2(x)\Vert\pi_{\mathcal{F}_i(x)}h\Vert^2\,d\mu(x)+B\int_{\sigma^c}w_i^2(x)\Vert\pi_{\mathcal{G}_i(x)}h\Vert^2\,d\mu(x)\\
&\quad\geq \int_{\sigma}\int_{X_i}\vert\langle \pi_{\mathcal{F}_i(x)}h, v_i(x)F_i(y)\rangle\vert^2\,d\mu(y)\,d\mu(x)\\
&\qquad +\int_{\sigma^c}\int_{X_i}\vert\langle \pi_{\mathcal{G}_i(x)}h, w_i(x)G_i(y)\rangle\vert^2\,d\mu(y)\,d\mu(x)\\
&\quad\geq C\Vert h\Vert^2,
\end{align*}
where $B=\max\{B_{F},B_{G}\}$. Thus, we obtain $(\mathcal{F}_i, v_i)_{i\in[m]}$ and $(\mathcal{G}_i, w_i)_{i\in[m]}$ are c-fusion woven for $H$ with universal frame bounds $\frac{C}{B}$ and $\frac{D}{A}$.

$(II)\Rightarrow(I)$. Suppose that $(\mathcal{F}_i, v_i)_{i\in[m]}$ and $(\mathcal{G}_i, w_i)_{i\in[m]}$ are c-fusion woven for $H$ with universal frame bounds $C$ and $D$. Now, we can write for each $h\in H$,
\begin{align*}
&\int_{\sigma}\int_{X_i}\vert\langle h, v_i(x)F_i(y)\rangle\vert^2\,d\mu(y)\,d\mu(x)
+\int_{\sigma^c}\int_{X_i}\vert\langle h, w_i(x)G_i(y)\rangle\vert^2\,d\mu(y)\,d\mu(x)\\
&\quad=\int_{\sigma}\int_{X_i}\vert\langle \pi_{\mathcal{F}_i(x)}h, v_i(x)F_i(y)\rangle\vert^2\,d\mu(y)\,d\mu(x)\\
&\qquad+\int_{\sigma^c}\int_{X_i}\vert\langle \pi_{\mathcal{G}_i(x)}h, w_i(x)G_i(y)\rangle\vert^2\,d\mu(y)\,d\mu(x)\\
&\quad \geq \int_{\sigma}A_{F_i}\Vert v_i(x)\pi_{\mathcal{F}_i(x)}h\Vert^2\,d\mu(x)+\int_{\sigma^c}A_{G_i}\Vert w_i(x)\pi_{\mathcal{G}_i(x)}h\Vert^2\,d\mu(x)\\
&\quad\geq A\Big(\int_{\sigma}v_i^2(x)\Vert\pi_{\mathcal{F}_i(x)}h\Vert^2\,d\mu(x)+\int_{\sigma^c}w_i^2(x)\Vert\pi_{\mathcal{G}_i(x)}h\Vert^2\,d\mu(x)\Big)\\
&\quad\geq AC\Vert h\Vert^2.
\end{align*}
Also, we get
\begin{small}
\begin{align*}
\int_{\sigma}\int_{X_i}\vert\langle h, v_i(x)F_i(y)\rangle\vert^2\,d\mu(y)\,d\mu(x)
+\int_{\sigma^c}\int_{X_i}\vert\langle h, w_i(x)G_i(y)\rangle\vert^2\,d\mu(y)\,d\mu(x)
\leq BD\Vert h\Vert^2.
\end{align*}
\end{small}
So, $\{\mathfrak{F}_i.\mathfrak{v}_i\}_{i\in[m]}$ and $\{\mathfrak{G}_i.\mathfrak{w}_i\}_{i\in[m]}$  are c-woven for $H$  with universal bounds $AC$ and $BD$.
\end{proof}
\begin{example}
Let $X=\Bbb R^2$, $\mu$ be the Lebesgue measure, $H=L^2(\Bbb R)$ and
$$X_1=\{(a,b): \ b\geq0\},\quad X_2=\{(a,b): \ b\leq0\}.$$
Suppose that
\begin{align*}
F_i,G_i&:X_i\longrightarrow H,\\
F_i(a,b)=\mu(X_i)E_{a}T_{b} g,&\qquad G_i(a,b)=\mu(X_i)E_{a}T_{b}(\alpha g),
\end{align*}
where, $\alpha$ is a scaler such that $\vert\alpha\vert^2>1$, $i\in\{1,2\}$ and $g\in L^2(\Bbb R)\setminus\{0\}$. Let
\begin{align*}
v_i, w_i&:X\longrightarrow\Bbb R^+,\\
v_i(a,b)=&w_i(a,b)=\frac{1}{\mu(X_i)}.
\end{align*}
Therefore,  we have for all measurable subset $\sigma\in X$ and $f\in H$,
\begin{small}
\begin{align*}
&\int_{\sigma}\int_{X_i}\vert\langle f, v_i(x)F_i(a,b) \rangle\vert^2\,d\mu(a,b)\mu(x)+\int_{\sigma^c}\int_{X_i}\vert\langle f, w_i(x)G_i(a,b) \rangle\vert^2\,d\mu(a,b)\mu(x)\\
&\quad=\int_{\sigma}\int_{X_i}\vert\langle f, E_aT_b g \rangle\vert^2\,d\mu(a,b)\mu(x)
+\vert\alpha\vert^2\int_{\sigma^c}\int_{X_i}\vert\langle f, E_aT_bg \rangle\vert^2\,d\mu(a,b)\mu(x)\\
&\quad\geq\int_{X}\int_{X_i}\vert\langle f, E_aT_bg \rangle\vert^2\,d\mu(a,b)\mu(x),
\end{align*}
\end{small}
where $i\in\{1,2\}$. Hence, by \cite{ch} and notaitions of Theorem \ref{tfF}, $\{\mathfrak{F}_i.\mathfrak{v}_i\}_{i\in\{1,2\}}$ and $\{\mathfrak{G}_i.\mathfrak{w}_i\}_{i\in\{1,2\}}$ are c-woven for $H$.
Thus, by Theorem \ref{tfF}, $(\mathcal{F}_i, v_i)_{i\in\{1,2\}}$ and $(\mathcal{G}_i, w_i)_{i\in\{1,2\}}$ are c-fusion woven for $H$.
\end{example}
\begin{theorem}\label{notw}
Let $(F_i, v_i)_{i\in[m]}$ be a family of c-fusion frame for $H$ with respect to a $\sigma$-finite measure $\mu$. Suppose that for any partition collection of disjoint finite sets $\{\tau_i\}_{i\in[m]}$ of $X$ and for any $\varepsilon>0$ there exists a partition $\{\sigma_i\}_{i\in[m]}$ of the set $X\setminus\bigcup_{i\in[m]}\tau_i$ such that $(F_i(x), v_i(x))_{x\in(\sigma_i\cup\tau_i), i\in[m]}$ has a lower c-fusion frame bound less than $\varepsilon$. Then
$(F_i(x), v_i(x))_{x\in X, i\in[m]}$ is not a c-fusion woven.
\end{theorem}
\begin{proof}
Since $(X,\mu)$ is a $\sigma$-finite measure space, then  $X=\cup_{i\in\Bbb N}X_i$, where $X_i$ are disjoint measurable sets and $\mu(X_i)<\infty$ for all $i\in\Bbb N$. Assume that $\tau_{i1}=\emptyset$ for all $i\in[m]$ and $\varepsilon=1$. Then, there exists a partition $\{\sigma_{i1}\}_{i\in[m]}$ of $X$ such that  $(F_i(x), v_i(x))_{x\in(\sigma_{i1}\cup\tau_{i1}), i\in[m]}$ has a lower bound (also, optimal lower bound) less that $1$. Thus, there is a vector $h_1\in H$ with $\Vert h_1\Vert=1$ such that
$$\sum_{i\in[m]}\int_{\sigma_{i1}\cup\tau_{i1}}v_{i}^2(x)\Vert\pi_{F_i(x)}h_1\Vert^2\,d\mu<1.$$
Since
$$\sum_{i\in[m]}\int_{X}v_{i}^2(x)\Vert\pi_{F_i(x)}h_1\Vert^2\,d\mu<\infty,$$
so, there is a $k_1\in\Bbb N$ such that
$$\sum_{i\in[m]}\int_{\Bbb K_1}v_{i}^2(x)\Vert\pi_{F_i(x)}h_1\Vert^2\,d\mu<1,$$
where, $\Bbb K_1=\cup_{i\geq k_1+1}X_i$.

Continuing this way, for $\varepsilon=\frac{1}{n}$ and a partition $\{\tau_{in}\}_{i\in[m]}$ of $X_1\cup\cdots\cup X_{k_n-1}$ such that
$$\tau_{in}=\tau_{i(n-1)}\cup\big(\sigma_{i(n-1)}\cap(X_1\cup\cdots\cup X_{k_n-1})\big)$$
for all $i\in[m]$, there exists a partition $\{\sigma_{in}\}_{i\in[m]}$ of $X\setminus(X_1\cup\cdots\cup X_{k_n-1})$ such that $(F_i(x), v_i(x))_{x\in(\sigma_{in}\cup\tau_{in}), i\in[m]}$ has a lower bound less than $\frac{1}{n}$. Therefore, there is a $h_{n}\in H$ and $k_n\in\Bbb N$ such that $\Vert h_n\Vert=1$, $k_n>k_{n-1}$ and
$$\sum_{i\in[m]}\int_{\Bbb K_n}v_{i}^2(x)\Vert\pi_{F_i(x)}h_n\Vert^2\,d\mu<\frac{1}{n},$$
where, $\Bbb K_n=\cup_{i\geq k_n+1}X_i$. Choose a partition $\{\varsigma_i\}_{i\in[m]}$ of $\Bbb J$, where $\varsigma_i:=\cup_{j\in\Bbb N}\{\tau_{ij}\}=\tau_{i(n+1)}\cup(\varsigma_{i}\cap X\setminus(X_1\cup\cdots\cup X_{n}))$. Assume that $(F_i(x), v_i(x))_{x\in\varsigma_i, i\in[m]}$ is a c-fusion frame for $H$ with the optimal lower frame bound $A$. Then, by the Archimedean Property, there exists a $n\in\Bbb N$ such that $r>\frac{2}{A}$. Now, there exists a $h_r\in H$ with $\Vert h_r\Vert=1$ such that
\begin{small}
\begin{align*}
&\sum_{i\in[m]}\int_{\varsigma_i}v_{i}^2(x)\Vert\pi_{F_i(x)}h_r\Vert^2\,d\mu\\
&\quad=\sum_{i\in[m]}\int_{\tau_{i(r+1)}}v_{i}^2(x)\Vert\pi_{F_i(x)}h_r\Vert^2\,d\mu+
\sum_{i\in[m]}\int_{\varsigma_{i}\cap X\setminus(X_1\cup\cdots\cup X_{r})}v_{i}^2(x)\Vert\pi_{F_i(x)}h_r\Vert^2\,d\mu\\
&\quad\leq\sum_{i\in[m]}\int_{\tau_{ir}\cup\sigma_{ir}}v_{i}^2(x)\Vert\pi_{F_i(x)}h_r\Vert^2\,d\mu+
\sum_{i\in[m]}\int_{\cup_{k\geq r+1}X_k}v_{i}^2(x)\Vert\pi_{F_i(x)}h_r\Vert^2\,d\mu\\
&\quad <\frac{1}{r}+\frac{1}{r}\\
&\quad <A\Vert h_r\Vert^2,
\end{align*}
\end{small}
and this is a contradiction with the lower bound of $A$.
\end{proof}
\begin{corollary}
Let $(F_i(x), v_i(x))_{x\in X, i\in[m]}$ be a c-fusion woven for $H$ with respect to a $\sigma$-finite measure $\mu$. Then there exists a collection of disjoint measurable subsets $\{\tau_i\}_{i\in[m]}$ of $X$ and $A>0$ such that for any partition $\{\sigma_i\}_{i\in[m]}$ of the set $X\setminus\cup_{i\in[m]}\tau_i$, the family $(F_i(x), v_i(x))_{x\in(\tau_i\cup\sigma_i), i\in[m]}$ is a c-fusion frame for $H$ with the lower frame bound $A$.
\end{corollary}
\begin{corollary}\label{opt}
Suppose that $(F,v)$ and $(G,w)$ are c-fusion frames for $H$ with optimal upper frame bounds $B_1$ and $B_2$, respectively, and they are c-fusion woven for $H$. Then, $B_1+B_2$ is not an optimal upper frame bound for the c-fusion woven.
\end{corollary}
\begin{proof}
Let $\varepsilon>0$. Assume that $B_1+B_2$ is an optimal upper frame bound  for the c-fusion woven. So, there exists $\sigma\subset X$ such that
$$\sup_{\Vert h\Vert=1}\Big(\int_{\sigma}v^2(x)\Vert\pi_{F(x)}h\Vert^2\,d\mu+\int_{\sigma^c}w^2(x)\Vert\pi_{G(x)}h\Vert^2\,d\mu\Big)=B_1+B_2.$$
Therefore, there is a $h_1\in H$ and $\Vert h_1\Vert=1$ such that
$$\int_{\sigma}v^2(x)\Vert\pi_{F(x)}h_1\Vert^2\,d\mu+\int_{\sigma^c}w^2(x)\Vert\pi_{G(x)}h_1\Vert^2\,d\mu\geq B_1+B_2-\varepsilon.$$
Thus, by the hypothesis,
$$\int_{X\setminus\sigma}v^2(x)\Vert\pi_{F(x)}h_1\Vert^2\,d\mu+\int_{X\setminus\sigma^c}w^2(x)\Vert\pi_{G(x)}h_1\Vert^2\,d\mu\leq\varepsilon.$$
Now, by Theorem \ref{notw}, we conclude that $(F,v)$ and $(G,w)$ are not c-fusion woven and this is a contradiction.
\end{proof}
\begin{theorem}
Let $(F_i(x), v_i(x))_{x\in X}$ be a c-fusion frame for $H$ with frame bounds $A_i$ and $B_i$ for each $i\in[m]$. Suppose that there exists $N>0$ so that for all $h\in H$, $i\neq k\in[m]$ and all measurable subset $Y\subset X$,
\begin{align*}
&\int_{Y}\Vert(v_i(y)\pi_{F_i(y)}-v_k(y)\pi_{F_k(y)})h\Vert^2\,d\mu(y)\\
&\qquad\leq N\min\Big\lbrace\int_{Y}\Vert v_i(y)\pi_{F_i(y)}h\Vert^2\,d\mu(y),\int_{Y}\Vert v_k(y)\pi_{F_k(y)}h\Vert^2\,d\mu(y)\Big\rbrace.
\end{align*}
Then, the family $(F_i(x), v_i(x))_{x\in X, i\in[m]}$ is a c-fusion frame with universal bounds
$$\frac{A}{(m-1)(N+1)+1}\quad \mbox{and}\quad B,$$
where, $A:=\sum_{i\in[m]}A_i$ and $B=\sum_{i\in[m]}B_i$.
\end{theorem}
\begin{proof}
Assume that $\{\sigma_i\}_{i\in[m]}$ is a partition of $X$ and $h\in H$. Thus,
\begin{small}
\begin{align*}
\sum_{i\in[m]}A_i\Vert h\Vert^2&\leq\sum_{i\in[m]}\int_{X}v_i^2(x)\Vert\pi_{F_i(x)}h\Vert^2\,d\mu\\
&=\sum_{i\in[m]}\sum_{k\in[m]}\int_{\sigma_k}v_i^2(x)\Vert\pi_{F_i(x)}h\Vert^2\,d\mu\\
&\leq\sum_{i\in[m]}\bigg(\int_{\sigma_i}v_i^2(x)\Vert\pi_{F_i(x)}h\Vert^2\,d\mu\\
&\quad+\sum_{\substack{k\in[m] \\ k\neq i}}\int_{\sigma_k}\Big\lbrace\Vert v_i^2(x)\pi_{F_i(x)}h-v_k^2(x)\pi_{F_k(x)}h\Vert^2\,d\mu
+v_k^2(x)\Vert\pi_{F_k(x)}h\Vert^2\,d\mu\Big\rbrace\bigg)\\
&\leq\sum_{i\in[m]}\Big(\int_{\sigma_i}v_i^2(x)\Vert\pi_{F_i(x)}h\Vert^2\,d\mu+
\sum_{\substack{k\in[m] \\ k\neq i}}\int_{\sigma_k}(N+1)v_k^2(x)\Vert\pi_{F_k(x)}h\Vert^2\,d\mu\Big)\\
&=\lbrace(m-1)(N+1)+1\rbrace\sum_{i\in[m]}\int_{\sigma_i}v_i^2(x)\Vert\pi_{F_i(x)}h\Vert^2\,d\mu.
\end{align*}
\end{small}
Hence, by Proposition \ref{t1}, the proof is completed.
\end{proof}
%%%%%%%%%%%%%%%%%%%%%%%%%%%%%%%%%%
\section{Perturbation For C-Fusion Woven}
\begin{theorem}\label{per}
For every $i\in[m]$, let $(F_i(x), v_i(x))_{x\in X}$ be a c-fusion frame for $H$ with frame bounds $A_i$ and $B_i$.  Suppose that there exist non-negative scalers $\lambda_i$, $\eta_i$ and $\gamma_i$ such that for some fixed $n\in[m]$,
$$A:=A_n-\sum_{i\in[m]\setminus\{n\}}(\lambda_i+\eta_i\sqrt{B_n}+\gamma_i\sqrt{B_i})(\sqrt{B_n}+\sqrt{B_i})>0$$
and
%\begin{small}
\begin{align*}
\Big\Vert\int_{X}(v_nf-v_if)\,d\mu\Big\Vert
\leq\eta_{i}\Big\Vert\int_{X}v_nf\,d\mu\Big\Vert
+\gamma_i\Big\Vert\int_{X}v_if\,d\mu\Big\Vert
+\lambda_i\Vert f\Vert
\end{align*}
%\end{small}
for every $f\in\mathscr{L}^2(X,F)$ and $i\in[m]\setminus\{n\}$. Then the family $(F_i(x), v_i(x))_{x\in X, i\in[m]}$ is a g-fusion woven for $H$ with universal frame bounds $A$ and $\sum_{i\in[m]}B_i$.
\end{theorem}
\begin{proof}
By Proposition \ref{t1}, the upper bound is obvious. For the lower bound, assume that $T_{F_i}$ is the synthesis operator of $(F_i, v_i)$, $f\in\mathscr{L}^2(X,F)$ and $i\in[m]\setminus\{n\}$. Therefore,
\begin{align*}
\Vert(T_{F_n}-T_{F_i})f\Vert&=\Vert\int_{X}(v_nf-v_if)\,d\mu\Vert\\
&\leq\eta_{i}\Big\Vert\int_{X}v_nf\,d\mu\Big\Vert
+\gamma_i\Big\Vert\int_{X}v_nf\,d\mu\Big\Vert
+\lambda_i\Vert f\Vert\\
&\leq(\eta_{i}\sqrt{B_n}+\gamma_i\sqrt{B_i}+\lambda_i)\Vert f\Vert.
\end{align*}
Thus,
\begin{align}\label{2}
\Vert T_{F_n}-T_{F_i}\Vert\leq(\eta_{i}\sqrt{B_n}+\gamma_i\sqrt{B_i}+\lambda_i).
\end{align}
For each $i\in[m]$ and $\sigma\subset X$, we define
\begin{align*}
T_{i}^{(\sigma)}&:\mathscr{L}^2(X,F) \to H,\\
T_{i}^{(\sigma)}&(f)=\int_{\sigma}v_if\,d\mu.
\end{align*}
Then, $T_{i}^{(\sigma)}=T_{F_i}.\chi_{\sigma}$. So, $T_{i}^{(\sigma)}$ is bounded and $\Vert T_{i}^{(\sigma)}\Vert\leq\sqrt{B_i}$ for all $i\in[m]$. Similarly with \eqref{2}, we can get for each $i\in[m]\setminus\{n\}$,
\begin{align}
\Vert T_{n}^{(\sigma)}-T_{i}^{(\sigma)}\Vert\leq(\eta_{i}\sqrt{B_n}+\gamma_i\sqrt{B_i}+\lambda_i).
\end{align}
We compute for every $h\in H$ and $i\in[m]\setminus\{n\}$,
\begin{align*}
\big\Vert&\big(T_n^{(\sigma)}(T_n^{(\sigma)})^*-T_i^{(\sigma)}(T_i^{(\sigma)})^*\big)h\big\Vert\\
&\leq\big\Vert\big(T_n^{(\sigma)}(T_n^{(\sigma)})^*-T_n^{(\sigma)}(T_i^{(\sigma)})^*\big)h\big\Vert+\big\Vert\big(T_n^{(\sigma)}(T_i^{(\sigma)})^*-T_i^{(\sigma)}(T_i^{(\sigma)})^*\big)h\big\Vert\\
&\leq\big\Vert T_n^{(\sigma)}\big\Vert\big\Vert \big((T_n^{(\sigma)})^*-(T_i^{(\sigma)})^*\big)h\big\Vert
+\big\Vert T_i^{(\sigma)}\big\Vert\big\Vert \big(T_n^{(\sigma)}-T_i^{(\sigma)}\big)h\big\Vert\\
&\leq(\eta_{i}\sqrt{B_n}+\gamma_i\sqrt{B_i}+\lambda_i)(\sqrt{B_n}+\sqrt{B_i})\Vert h\Vert.
\end{align*}
Now, suppose that $\{\sigma_i\}_{i\in[m]}$ is a partition of $X$ and $T_F$ is the synthesis operator associated with the c-fusion Bessel sequence $(F_i(x), v_i(x))_{x\in X, i\in[m]}$. Hence, we obtain for any $h\in H$,
\begin{align*}
\Vert T^*_F h\Vert^2&=\vert\langle h, T_FT^*_F h\rangle\vert\\
&=\Big\vert\Big\langle h,\sum_{i\in[m]}\int_{\sigma_i}v_i^2(x)\pi_{F_i(x)}h\,d\mu(x)\Big\rangle\Big\vert\\
&=\Big\vert\Big\langle h,
\sum_{i\in[m]\setminus\{n\}}\int_{\sigma_i}v_i^2(x)\pi_{F_i(x)}h\,d\mu(x)+
\sum_{i\in[m]}\int_{\sigma_i}v_n^2(x)\pi_{F_n(x)}h\,d\mu(x)\\
&\quad-\sum_{i\in[m]\setminus\{n\}}\int_{\sigma_i}v_n^2(x)\pi_{F_n(x)}h\,d\mu(x)\Big\rangle\Big\vert\\
&=\Big\vert\Big\langle h,
\sum_{i\in[m]}\int_{\sigma_i}v_n^2(x)\pi_{F_n(x)}h\,d\mu(x)\\
&\quad-\sum_{i\in[m]\setminus\{n\}}\int_{\sigma_i}\Big(v_n^2(x)\pi_{F_n(x)}-v_i^2(x)\pi_{F_i(x)}\Big)h\,d\mu(x)\Big\rangle\Big\vert\\
&\geq\Big\vert\Big\langle h,\sum_{i\in[m]}\int_{\sigma_i}v_n^2(x)\pi_{F_n(x)}h\,d\mu(x)\Big\rangle\Big\vert\\
&\quad -\sum_{i\in[m]\setminus\{n\}}\Big\vert\Big\langle h,\int_{\sigma_i}\Big(v_n^2(x)\pi_{F_n(x)}-v_i^2(x)\pi_{F_i(x)}\Big)h\,d\mu(x)\Big\rangle\Big\vert\\
&\geq\vert\langle h, T_{F_n}T^*_{F_n} h\rangle\vert-\sum_{i\in[m]\setminus\{n\}}\Vert h\Vert \big\Vert \big(T_n^{(\sigma_i)}(T_n^{(\sigma_i)})^*-T_i^{(\sigma_i )}(T_i^{(\sigma_i)})^*\big)h\big\Vert\\
&\geq A_n\Vert h\Vert^2-\sum_{i\in[m]\setminus\{n\}}\Vert h\Vert(\eta_{i}\sqrt{B_n}+\gamma_i\sqrt{B_i}+\lambda_i)(\sqrt{B_n}+\sqrt{B_i})\Vert h\Vert\\
&=\Big(A_n-\sum_{i\in[m]\setminus\{n\}}(\eta_{i}\sqrt{B_n}+\gamma_i\sqrt{B_i}+\lambda_i)(\sqrt{B_n}+\sqrt{B_i})\Big)\Vert h\Vert^2\\
&=A\Vert h\Vert^2.
\end{align*}
This completes the proof.
\end{proof}
With similar proof of Theorem \ref{per}, we can show the following result  when the index $n$  is not fixed.
\begin{corollary}
For each $i\in[m]$, let $(F_i(x), v_i(x))_{x\in X}$ be a c-fusion frame for $H$ with frame bounds $A_i$ and $B_i$.  Suppose that there exist non-negative scalers $\lambda_i$, $\eta_i$, $\gamma_i$ and $n\in[m]$ so that,
$$A:=A_1-\sum_{i\in[m-1]}(\lambda_i+\eta_i\sqrt{B_i}+\gamma_i\sqrt{B_{i+1}})(\sqrt{B_i}+\sqrt{B_{i+1}})>0$$
and
%\begin{small}
\begin{align*}
\Big\Vert\int_{X}(v_if-v_{i+1}f)\,d\mu\Big\Vert
\leq\eta_{i}\Big\Vert\int_{X}v_if\,d\mu\Big\Vert
+\gamma_i\Big\Vert\int_{X}v_{i+1}f\,d\mu\Big\Vert
+\lambda_i\Vert f\Vert
\end{align*}
%\end{small}
for every $f\in\mathscr{L}^2(X,F)$ and $i\in[m-1]$. Then the family $(F_i(x), v_i(x))_{x\in X, i\in[m]}$ is a g-fusion woven for $H$ with universal frame bounds $A$ and $\sum_{i\in[m]}B_i$.
\end{corollary}

%%%%%%%%%%%%%%%%%%%%%%%%%%%%%%%%%%%%%%%%

\end{document}